\magnification=\magstep1
\input amstex
\documentstyle{amsppt}
\NoBlackBoxes
\NoRunningHeads
\nologo
\define\qbinnk{\left[\matrix n \\ k \endmatrix \right]_q}
\topmatter
\title{Specializations of generalized Laguerre polynomials}
\endtitle
\author R. Simion\footnote
{This work was carried out in part during the first author's visits at the 
Mittag-Leffler Institute and the University of Qu\'ebec at Montr\'eal,
and with partial support through NSF grant DMS-9108749.}
and D. Stanton
\footnote{This work was supported by the 
Mittag-Leffler Institute and by NSF grant DMS-9001195.}
\endauthor
\address Department of Mathematics, George Washington University,
Washington, DC 20052.
\endaddress
\address School of Mathematics, University of Minnesota,
Minneapolis, MN 55455.
\endaddress
\dedicatory The second author dedicates this paper to the one and only Dick Askey, his mathematical father.
\enddedicatory
\abstract
Three specializations of a set of orthogonal polynomials with ``8 
different q's'' are given. The polynomials are identified as 
$q$-analogues of Laguerre polynomials, and the combinatorial 
interpretation of the moments give infinitely many new 
Mahonian statistics on permutations. 
\endabstract
\endtopmatter

\subheading{1. Introduction}

The Laguerre polynomials $L_n^{\alpha}(x)$ have been extensively studied, 
analytically \cite{E} and combinatorially \cite{F-S},\cite{V1}. 
There is also a 
classical set of $q$-Laguerre polynomials, \cite{M}. Recently, a set of 
orthogonal polynomials generalizing the Laguerre polynomials 
has been studied \cite{Si-St}. These polynomials in some sense have 
``8 different q's". Various specializations of them give 
orthogonal polynomials 
associated with many types of combinatorial objects. The 
purpose of this paper is to present the specializations which are true 
$q$-analogues of $L_n^0(x)$. By this we mean that the $n^{th}$ moments, 
instead of being $n!$, are basically $n!_q$.

We present here three specializations whose moments lead to new Mahonian 
statistics on permutations (Theorems 2, 3, and 4). 
In fact, infinitely many Mahonian statistics can be derived from 
those presented here.  Moreover, the theorems obtained from our 
specializations follow easily from classical analytic facts,
but are combinatorially non-trivial.

We shall use the terminology and notation found in \cite{G-R}, and let
$$
[n]_q=\frac{1-q^n}{1-q},\quad [n]_{r,s}= \frac{r^n-s^n}{r-s}.
$$

\subheading{2. The polynomials and their moments}

Any set of monic orthogonal polynomials satisfies the three term 
recurrence relation
$$
p_{n+1}(x) = (x-b_n)p_n(x)-\lambda_np_{n-1}(x).
\tag2.1a
$$
For the set of orthogonal polynomials with 8 different ``q's" considered in \cite{Si-St}, the coefficients are
$$
b_n=a[n+1]_{r,s}+b[n]_{t,u}, \quad \lambda_n=ab[n]_{p,q}[n]_{v,w}.
\tag2.1b
$$
We refer to the polynomials defined by (2.1) as the {\it octabasic 
Laguerre polynomials}.  They generalize the Laguerre polynomials 
and are the polynomials whose specializations we consider in this paper.

The fundamental combinatorial fact (Theorem 1) 
that we need here
concerns the moments for these polynomials. They are 
generating functions 
for permutations according to certain statistics.

For the definition
of these statistics, it is convenient to represent a
permutation $\sigma$ as a word $\sigma(1)\sigma(2)\cdots\sigma(n)$
consisting of increasing runs, separated by the
descents of the permutation. For example, the permutation
$\sigma=26|357|4|189$ has 4 runs separated by 3 descents,
and we write run$(\sigma)$ = 3.
The runs of length $2$ or more will be called {\it proper runs}
and those of length $1$ will be called {\it singleton runs}.

The elements $\sigma(i)$ of $\sigma$ fall into four classes:
the elements which begin proper runs (openers),
the elements which close proper runs (closers), the elements which
form singleton runs (singletons), and the elements which continue runs
(continuators). We shall abbreviate these classes of
elements ``op", ``clos", ``sing", and ``cont" respectively.
In the example, $op(\sigma) = \{3,2,1\}$,
$clos(\sigma) = \{6,7,9\}$, $sing(\sigma) =  \{4\}$, and
$cont(\sigma) =  \{5,8\}$.

\proclaim{Definition 1}
For $\sigma\in S_n$,  the statistics $lsg(\sigma)$ and
$rsg(\sigma)$ are defined by
$$
lsg(\sigma) = \sum_{i=1}^n lsg(i), \qquad rsg(\sigma) = \sum_{i=1}^n rsg(i),
$$
where $lsg(i)=$ the number of runs of $\sigma$ strictly
to the left of $i$ which contain elements smaller and greater than $i$,
and $rsg(i)=$ the number of runs of $\sigma$ strictly
to the right of $i$ which contain
elements smaller and greater than $i$.

We also define the $lsg$ and $rsg$ on the openers of $\sigma$
$$
lsg(op)(\sigma) = \sum_{i \in op(\sigma)} lsg(i),
\qquad  rsg(op)(\sigma) = \sum_{i \in op(\sigma)} lsg(i).
$$
The statistics  lsg  and   rsg  have analogous definitions on 
each of the remaining three classes of elements.
\endproclaim

For example, if  $\sigma=26|357|4|189$, then
$lsg(7) = 0$, $rsg(7) = 1$, $lsg(op)(\sigma) = 0+1+0=1$,
$rsg(op)(\sigma) = 1+1+0 = 2$,
$lsg(clos)(\sigma) = 0$,
$rsg(clos)(\sigma) = 2+1+0 = 3$, etc.

\proclaim{Theorem 1}
The $n^{th}$ moment $\mu_n$ for the octabasic Laguerre polynomials is
$$
\align
\mu_n=\sum_{\sigma\in S_n}&
r^{{\text{lsg(sing)}}(\sigma)}s^{{\text{rsg(sing)}}(\sigma)}
t^{{\text{lsg(cont)}}(\sigma)}u^{{\text{rsg(cont)}}(\sigma)}
p^{{\text{lsg(op)}}(\sigma)}q^{{\text{rsg(op)}}(\sigma)}\\
&v^{{\text{lsg(clos)}}(\sigma)}w^{{\text{rsg(clos)}}(\sigma)}
a^{{\text{run}}(\sigma)}b^{n-{\text{run}}(\sigma)}.
\endalign
$$
\endproclaim

\demo{Proof}
The Viennot theory \cite{V1}, 
\cite{V2} gives $\mu_n$ as a generating function for
Motzkin paths from $(0,0)$ to $(n,0)$. From (2.1b) the paths have 4 types
of edges (or steps):
\roster
\item northeast (NE) edges starting at level $k$ with weight in $a[k+1]_{p,q}$,
\item southeast (SE) edges starting at level $k$ with weight in $b[k]_{v,w}$,
\item east (E solid) edges starting at level $k$ with weight in $a[k+1]_{r,s}$,
\item east (E dotted) edges starting at level $k$ with weight in $b[k]_{t,u}$,
\endroster
where ``weight in $\alpha[m]_{c,d}$'' means that the weight
of the edge is one of the monomials
$\alpha c^{m-1}, \alpha c^{m-2}d, \dots, \alpha d^{m-1}$ which
appear in  $\alpha[m]_{c,d}$.

The weight of a Motzkin path is then the product of
the weights of its edges.
An example of such a weighted path is given in Figure 1.

In the statement of Theorem 1 we claim that $\mu_n$ is
the sum over all permutations in $S_n$ of monomial
weights defined in terms of our permutation statistics.
To prove the theorem we construct a bijection between
the weighted Motzkin paths of length $n$
and permutations in $S_n$,
so that the weight of a path is equal to the weight of
its corresponding  permutation.

Given a weighted Motzkin path of length $n$,
its corresponding permutation will be constructed in
$n$ stages.
We begin at the origin and with the empty permutation.
We traverse the path from left to right,
and the $i^{th}$ step will determine where to insert $i$ in the
current (partial) permutation $\sigma_{i-1} \in S_{i-1}$.
The end result will be a permutation $\sigma_n = \sigma \in S_n$.
Depending on the type of the $i^{th}$ step of the path,
$i$ will belong to one or another of the four classes of
elements of $\sigma$:
\roster
\item   NE:  $i \in \ op(\sigma)$,
\item   SE:  $i \in \ clos(\sigma)$,
\item   E solid:  $i \in \ sing(\sigma)$,
\item   E dotted: $i \in \ cont(\sigma)$.
\endroster

The weight $\alpha c^j d^k$ of the $i^{th}$ step of the path determines,
as described below, the exact position in
$\sigma_{i-1}$ where we insert $i$ as a point in the
appropriate class of elements.

A run in $\sigma_{i-1}$ will be called an {\it active run}
if its maximum is a (future)  opener or continuator in $\sigma$.
Notice that the first step of the path starts at level $0$
and that $\sigma_0$, being the empty permutation, has no active runs.
Inductively, assume that the level $h$ at which the
$i^{th}$ step starts is equal to the number of
active runs in $\sigma_{i-1}$,
and recall the relation between the weight $\alpha c^j d^k$
of a step and its starting level $h$:
for NE and E solid steps we have $j + k = h$,
while for other steps we have $j + k = h -1$.

If the $i^{th}$ step is  E dotted or SE,
then the element $i$ is adjoined to the
$(j+1)^{st}$ active run of $\sigma_{i-1}$,
as a continuator or  closer, respectively.
This is well-defined since, as discussed above,
$j+k$ is one unit less than the number of active
runs in $\sigma_{i-1}$.
It follows as well that the new partial permutation,
$\sigma_i$ will have as many active runs as the starting
level of the $(i+1)^{st}$ step of the path.

If the $i^{th}$ step is NE or E solid, then
$j+k$ is equal to the number of active runs in $\sigma_{i-1}$
and we insert $i$ in the leftmost position possible
such that it will have
$j$ of the {\it active} runs of $\sigma_{i-1}$  strictly
to its left and $k$ of the {\it active} runs of $\sigma_{i-1}$
strictly to its right.
The position where $i$ is inserted in this case
ensures that $i$ will be the initial element of a run,
and it is again true that $\sigma_i$ will have as many
active runs as the starting level of the $(i+1)^{st}$ step
of the path.

Note that in the final permutation $\sigma$
the values  of $lsg(i)$ and $rsg(i)$ are completely determined
by the runs which were active in $\sigma_{i-1}$ and the
position, relative to these runs, where $i$ was inserted.
So, if the $i^{th}$ step has weight $\alpha c^j d^k$,
then $lsg(i) = j$ because each of the $j$ open runs
of $\sigma_{i-1}$ which remain to the left of $i$ upon its
insertion will eventually be extended by at least one element
greater than $i$.  Similarly, $rsg(i) = k$.
Finally, $\alpha = a$ in the weight of the $i^{th}$ step corresponds
precisely to $i$ being the initial (possibly the only) element
of a run.
Consequently, the weight of a Motzkin path is equal
to our intended weight for the corresponding permutation $\sigma$.

It remains to verify that this correspondence is bijective.
We claim that, given $\sigma \in S_n$,
we can reconstruct its associated
Motzkin path, step by step, beginning at the right end of the path, $(n,0)$,
since each stage of our construction is reversible.

Using our rules (1)-(4), each permutation in $S_n$ will produce a
Motzkin path of length $n$ (not yet weighted),
since $|op(\sigma)| = |clos(\sigma)|$ and each closer is larger
than its corresponding opener.
We must check that the weight $\alpha c^{lsg(i)} d^{rsg(i)}$
(where $\alpha = a$ if $i$ is the initial element of a run,
and $\alpha = b$ otherwise) is a valid weight for
the $i^{th}$ step of the path.
This follows immediately  from the equality
$lsg(i)+rsg(i) = $ the level of the left endpoint $e_{i+1}$ of the partial
path reconstructed from the values $n, n-1, \dots, i+1$.
The equality holds for $i=n$ since $lsg(n) = rsg(n) =0$
in all permutations of $S_n$, and the left endpoint
of the one-point path consisting just of $(n,0) = e_{n+1}$ is at level $0$.
Suppose the equality holds for $i+1$ and we will prove it for $i$.
Observe that the level of $e_{i+1}$ is equal to the number of
SE steps minus the number of NE steps in the partial path
from $e_{i+1}$ to $(n,0)$.  That is,
the level of $e_{i+1}$ is equal to the number of
proper runs in $\sigma$ whose maximum is larger than $i$,
minus the number of proper runs in $\sigma$ whose minimum is larger than $i$.
But this is equal in turn with the number of proper runs in $\sigma$
with maximum larger than $i$ and minimum smaller than $i$,
i.e., it is equal to $lsg(i)+rsg(i)$.
It now becomes clear that our map from weighted Motzkin paths
to permutations is indeed invertible.
\qed
\enddemo

Figure 1 shows the weighted Motzkin path which corresponds with
the permutation $\sigma= 10|8 9 \thinspace 11| 1 3 7| 5| 4 6| 2$.
We have also included a binary tree representation of the
permutation, deeming it of possible interest to the readers
familiar with \cite{V1}. 
The bijection constructed in the proof above is related to
Viennot's correspondence between Motzkin paths and permutations
\cite{V1}.  As an intermediate step in Viennot's correspondence,
Motzkin paths and permutations are encoded by
increasingly labeled  binary trees.

We will be concerned with specializations of the 8 $q$'s under which the 
moments in Theorem 1 become multiples of $n!_q = [n]_q [n-1]_q \cdots [1]_q$. 
It is clear that the parameter $b$ can be rescaled to 1. 
Also, in considering specializations, we can take advantage of the 
property  ---  obvious from (2.1) --- that the moments are fixed 
under the interchange of $\{ r,s \}$, $\{ t,u \}$,  $\{ p, q\}$,
and   $\{ v, w\}$, and also fixed if $p$ and $q$ are interchanged with 
$v$ and $w$.

\subheading{3. The specializations}
 
In this section we state three different specializations of the polynomials 
defined by (2.1).  The polynomials which arise are monic little $q$-Jacobi,
sums of two little $q$-Jacobi, and classical $q$-Laguerre.
Each of these three cases will have moments which are basically 
$n!_q$.

First we choose the parameters so that the polynomials 
coincide with the monic form of the 
little $q$-Jacobi polynomials \cite{G-R, p. 166}, 
$p_n(xq(1-q);q^\alpha,0;q)$, which have 
$$
b_n=q^{n-1}[n+1+\alpha]_q+q^{n+\alpha-1}[n]_q,\quad \lambda_n=q^{2n-3+\alpha}[n]_q[n+\alpha]_q.
\tag3.1
$$
The appropriate specialization occurs only for 
$\alpha=0$, $p_n(xq(1-q);1,0;q)$, 
and is $r=t=p=v=q^2$, $s=u=q=w$, $a=1/q$, $b=1$. 
The explicit formula for the polynomials, which is 
just the definition of the little $q$-Jacobi polynomials, is
$$
p_n(x)=\sum_{k=0}^n \qbinnk [n]_q\cdots[n-k+1]_q (-1)^k 
x^{n-k} q^{\binom{k-1}{2}-1}.
\tag3.2 
$$

The measure for $p_n(x;q^{\alpha},q^{\beta};q)$ is purely discrete, 
with masses of 
$$
\frac{q^{(\alpha+1)i}(q^{\beta+1};q)_i 
(q^{\alpha+1};q)_\infty}{(q;q)_i(q^{\alpha+\beta+2};q)_\infty}
$$ 
at $x= q^i$. An easy calculation shows that the moments are given by
$\mu_n= \frac{(q^{\alpha+1};q)_n}{(q^{\alpha+\beta+2};q)_n}$. 
Thus for $p_n(xq(1-q);1,0;q)$, we have $\mu_n= q^{-n} n!_q$. 
Based upon these remarks, Theorem 1 becomes the following theorem. 
An equivalent theorem has been given in \cite{deM-V, Prop. 5.2}.

\proclaim{Theorem 2}
For $\sigma\in S_n$, let 
$$
s(\sigma):=n-{\text{run}}(\sigma)+2lsg(\sigma)+rsg(\sigma).
$$
Then
$$
\sum_{\sigma\in S_n} q^{s(\sigma)} = n!_q.
$$
\endproclaim
\vskip 8pt

Moreover, we see from the symmetry of (2.1) 
with respect to the 4 pairs of ``q's" that 
Theorem 2 holds for 16 statistics related to $s(\sigma)$. 
These 16 are obtained by choosing the coefficients 1 and 2 for 
$lsg$ and $rsg$ independently for the four types of elements of $\sigma$.
This means for example that 
$$
\align
s'(\sigma) =  n-{\text{run}}(\sigma)&+{\text{lsg(sing)+2rsg(sing)+lsg(op)+2rsg(op)}}\\
&+{\text{2lsg(cont)+rsg(cont)+2lsg(clos)+rsg(clos)}}
\endalign
$$
also satisfies Theorem 2. We will see later that 
in fact there are infinitely many equidistributed 
statistics related to $s(\sigma)$.

For our second specialization we consider
$$
b_n=q^{n+1}[n+1]_q+q^{n-1}[n]_q,\quad \lambda_n=q^{2n-1}[n]_q[n]_q.
\tag3.3
$$
The appropriate values are $r=t=p=v=q^2$, 
$s=u=q=w$, $a=q$, $b=1$. The polynomials 
turn out to be a sum of two little $q$-Jacobi polynomials,
$$
\align
p_n(x) = n!_q q^{\binom{n}{2}}(-1)^n\Bigl[
&{}_{2}\phi_1\left(
\matrix
q^{-n}&0;&q,&xq(1-q)\\
&q&&
\endmatrix
\right)\\
-(1-q^n)
&{}_{2}\phi_1\left(
\matrix
q^{1-n}&0;&q,&xq(1-q)\\
&q^2&&
\endmatrix
\right)\Bigr],
\endalign
$$
or equivalently
$$
p_n(x)=x^n+\sum_{k=1}^n \qbinnk [n]_q\cdots[n-k+2]_q ([n-k]_q+q^n) 
(-1)^k x^{n-k} q^{\binom{k}{2}} +(-1)^nq^{\binom{n+1}{2}}n!_q.
\tag3.4
$$
We omit the proof of these formulas. It is a verification of the 
recurrence relation (2.1) from the recurrence relation for the little 
$q$-Jacobi polynomials.

Since these polynomials do not explicitly appear in the literature, 
we cannot compute the 
moments by quoting the relevant facts about their measure. Nonetheless, the moments and measure are easily determined.

\proclaim{Proposition 1}
The moments for the polynomials in (3.3) are $\mu_0=1$, and 
$\mu_n=q\thinspace n!_q, n>0$.
The measure is purely discrete, with masses of 
$q^{i}(q;q)_\infty/(q;q)_{i-1}$ at $q^{i-1}/(1-q)$, $i\ge1$, and a 
mass of $1-q$ at 0. 
\endproclaim

\demo{Proof}
The $q$-binomial theorem clearly shows that the total mass $\mu_0=(1-q)+q=1$.
It also implies
$$
\align
\mu_n=&\sum_{i=1}^\infty \frac{q^{in}q^i(q;q)_\infty}{(1-q)^n(q;q)_{i-1}}
+(1-q)\delta_{n,0}\\
=&\thinspace qn!_q+(1-q)\delta_{n,0}.
\endalign
$$
Thus the stated measure has the right moments. 
To show that the polynomials are orthogonal with respect to this measure, 
note that it is easy to check, from the explicit formula (3.4), 
that the moments annihilate $p_1 , p_2 , \dots$.  Hence, the linear 
functional defined by the measure annihilates $p_1 , p_2 , \dots$. 
Finally, the three terms recurrence now shows that the polynomials 
are indeed orthogonal. 
\qed
\enddemo

We then get a companion theorem to Theorem 2. 
As in the case of Theorem 2, we have 16 equivalent versions of the 
statistic $s(\sigma)$, by assigning coefficients 1 and 2  
to $lsg$ and $rsg$ independently on openers, continuators, closers and 
singletons.
 
\proclaim{Theorem 3}
For $\sigma\in S_n$, let 
$$
s(\sigma):={\text{run}}(\sigma)-1+2lsg(\sigma)+rsg(\sigma).
$$
Then
$$
\sum_{\sigma\in S_n} q^{s(\sigma)} = n!_q.
$$
\endproclaim
\vskip 8pt

We remark that Theorems 2 and 3 are valid for 
an infinite number of variations of the statistic $s(\sigma)$. 
It is easy to verify that for each $\sigma\in S_n$
$$
lsg(op)(\sigma)+rsg(op)(\sigma)=lsg(clos)(\sigma)+rsg(clos)(\sigma).
\tag3.5
$$
(In fact there is a specialization of $\{r, s, t, u, p, q, v, w\}$
 with one free parameter giving (3.5).) 
Therefore the value of $s(\sigma )$ remains the same if the coefficients 
$\{ 1, 2 \}$ are replaced on the openers with  $\{ 1+c , 2+c \}$,
and on the closers with $\{ 1-c , 2-c \}$.
This provides a variation of Theorems 2 and 3 for each choice of 
the real parameter $c$. 
For example, $c=1$ gives the unusual choice of coefficients
$\{2,3\}$ and $\{0,1\}$.

Our third choice for specialization is the set of the classical 
$q$-Laguerre polynomials $L_n^\alpha(x(1-q);q)$ 
\cite{M}, \cite{G-R, p. 194}, whose monic form has 
$$
b_n=q^{-2n-\alpha}[n]_q+q^{-2n-1-\alpha}[n+1+\alpha]_q,\quad 
\lambda_n=q^{1-4n-2\alpha}[n]_q[n+\alpha]_q.
\tag3.6
$$
The appropriate values are $r=t=p=v=q^{-2}=b$, 
$s=u=q=w:=q^{-1}=a$ for
$L_n^0(x(1-q),q)$. 
The explicit formula for the monic form of $L_n^\alpha(x(1-q);q)$ is
$$
p_n(x)=\sum_{k=0}^n \qbinnk [n+\alpha]_q\cdots [n+\alpha-k+1]_q 
(-1)^kx^{n-k}q^{k(k-\alpha-2n)}.
\tag3.7
$$
Again a measure of these polynomials is explicitly known 
\cite{M, Th.1}, 
and the moments for $L_n^\alpha(x(1-q);q)$ can be found as 
$$
\mu_n=(q^{\alpha+1};q)_n q^{-n\alpha-\binom{n+1}{2}}/(1-q)^n.
$$
For $\alpha=0$ this is $q^{-\binom{n+1}{2}}n!_q$. However, the 
combinatorial 
version of this theorem is equivalent to Theorem 2, if $q$ is 
replaced with $1/q$.
Thus no new combinatorial theorem results. 

\subheading{4. The ``odd'' polynomials}

If $r=p$, $s=q$, $t=v$, and $u=w$, then the 
polynomials defined by (2.1) are the ``even" polynomials for the 
polynomials defined by (see \cite{C, p.41})
$$
b_n=0,\quad \lambda_{2n}=b[n]_{t,u},\quad \lambda_{2n+1}=a[n+1]_{r,s}.
$$
The ``odd" polynomials have the coefficients
$$
b_n=a[n+1]_{r,s}+b[n+1]_{r,s},\quad \lambda_n=ab[n+1]_{r,s}[n]_{t,u}.
\tag4.1
$$  
The moments for these ``odd" polynomials satisfy $\mu_n(odd)=\mu_{n+1}(even)/\mu_1(even)$. 
Since all of our 
specializations in \S3 satisfied $r=p$, $s=q$, $t=v$, and $u=w$, these 
``odd" polynomials also have moments which are multiples of $(n+1)!_q$. 
There is a version of Theorem 1 for the ``odd" polynomials which 
yields more statistics related to $s(\sigma)$. We do not state 
this combinatorial theorem here, rather in this section we state what 
these ``odd" polynomials are, give their moments, and state in Theorem 4 
what the
statistics related to $s(\sigma)$ are. 
Clearly the 
``odd" polynomials are analogues of the Laguerre polynomials $L_n^1(x)$.

We keep the parameters $r,s,t,u$.
This specialization gives the ``even'' and ``odd" 
polynomials a combinatorial 
interpretation as weighted versions 
of injective maps (see \cite{F-S}, $L_n^0(x)$ and $L_n^1(x)$). 
This family with ``4 q's" 
also contains other families of orthogonal polynomials of 
combinatorial interest which are discussed in \cite{Si-St}.

We list here the ``odd" polynomials for the three cases in \S3, and the 
respective moments. In each case the polynomials are monic forms 
of the given polynomials.

\roster
\item little $q$-Jacobi $p_n(xq(1-q);q,0;q)$, 
$\mu_n=q^{-n}\thinspace (n+1)!_q$,
\item little $q$-Jacobi $p_n(x(1-q);q,0;q)$, 
$\mu_n=(n+1)!_q$,
\item $q$-Laguerre $L_n^1(x(1-q);q)$,
$\mu_n=q^{-(n^2+3n)/2}\thinspace (n+1)!_q$.
\endroster

The combinatorial theorem that results is Theorem 4 below.  
To describe the suitable statistic $s(\sigma)$ we shall need an 
auxillary statistic, $n(\sigma )$, defined as follows. 
For a given permutation $\sigma \in S_n$, let $d$ be the largest 
element in the same run as 1, $d \ge 1$.  Now partition the elements of 
$\sigma$ into three classes:  elements to the left of 1, 
elements in the same run as 1, and those to the right of $d$. 
Suppose that over the portion of $\sigma$ to the right of $d$ 
no left-to-right minimum constitutes a singleton run. 
Put $n(\sigma ) = 0$.  Otherwise, let $c$ be the last 
(rightmost, smallest) singleton which is a left-to-right minimum 
on the portion of $\sigma$ to the right of $d$.  
Let $nleft(\sigma) \colon = \# \{ i \ \colon \ i < \sigma^{-1}(1), 
\ c < \sigma(i) < d \}$. 
In this case put $n(\sigma) = 2 (d-c) - nleft(\sigma)$. 
For example, $n(9\ |\ 1\ 5\ 7\ |\ 2\ 6\ |\ 4\ | \ 3 \ 8) = 0$, while for 
$\sigma = 7\ 12\ | \ 1\ 6\ 9\ | \ 3\ | \ 2\ 10\ 11\ |\ 5 \ | \ 4 \ 8$ 
we have $d=9$, $c=3$, $nleft(\sigma) = 1$ (from the element $7$), 
and we get $n(\sigma) = 2(9-3) -1 = 11$. 

We also need variations $lsg^*$ and $rsg^*$  on the statistics 
$lsg$ and $rsg$.   These differ from the original statistics in two 
respects.   For closers and singletons (the maxima of the runs), 
the run containing the element 1 is ignored 
in the calculation of $lsg^*$ and $rsg^*$.  For openers and continuators, 
the run containing 1 is always counted in $lsg^*$ (if it is to the left) 
or in $rsg^*$ (if it is to the right).  

For example, if 
$\sigma = 7\ 12\ | \ 1\ 6\ 9\ | \ 3\ | \ 2\ 10\ 11\ |\ 5 \ | \ 4 \ 8$,
then $lsg^*(5) = 1$, $lsg^*(10) = 2$, $lsg^*(\sigma) = 10$, 
$rsg^*(\sigma) = 8$.

\proclaim{Theorem 4}
For $\sigma\in S_n$, let 
$$
s(\sigma):={\text{run}}(\sigma)-1+2{\text{lsg*}}(\sigma)
+{\text{rsg*}}(\sigma)+n(\sigma).
$$
Then
$$
\sum_{\sigma\in S_n} q^{s(\sigma)} = n!_q.
$$ 
\endproclaim

Theorem 4 also holds if $run(\sigma)-1$ is replaced by $n-run(\sigma)$. 
Moreover, the role of closers and openers can be interchanged in Theorem 4, 
and there are also 16 variations, although complicated ones. 

A version of Theorem 4 holds for permutations in $S_{n+1}$ which 
satisfy the following condition:  $1$ and $n+1$ belong to the same run, 
and no left-to-right minimum to the right of $n+1$ constitutes a singleton 
run.  
There are $n!$ such 
permutations in $S_{n+1}$.

Finally, we remark that these specializations are the only ones we
have found for which the moments factor into an analogue of $n!$. They are 
also the only specializations which give the three sets of polynomials 
that were considered. A more extensive study of the specializations of 
(2.1) appears in \cite{Si-St}. 

\vskip1.5pc
{\bf References}
\item{[C]} T. Chihara, An Introduction to Orthogonal Polynomials,
 Mathematics and its applications, v. 13, 
Gordon and Breach, New York, 1978.
\item{[deM-V]} A. De Medicis and X. Viennot, Moments des $q$-polyn\^omes
de Laguerre et la bijection de Foata-Zeilberger, to appear.
\item{[E]} A. Erd\'elyi, Higher Transcendental Functions, 
McGraw-Hill, New York, 1953.
\item{[Fl]} P. Flajolet, Combinatorial aspects of continued fractions,
Discrete Math. 32 (1980), 126-161.
\item{[F-S]} D. Foata and V. Strehl, Combinatorics of Laguerre polynomials,
in Enumeration and Design, Waterloo Jubilee Conference, Academic Press, 
1984, 123-140.
\item{[G-R]} G. Gasper and M. Rahman, 
Basic Hypergeometric Series,  Encyclopedia of mathematics and its
applications, v. 35, Cambridge University Press, New York, 1990.
\item{[M]} D. Moak, The $q$-analogue of the Laguerre polynomials, J. Math.
Anal. Appl. 81 (1981), 20-47.
\item{[Si-St]} R. Simion and D. Stanton, Octabasic Laguerre polynomials 
and permutation statistics, in preparation.
\item{[V1]} G. Viennot, Une th\'eorie combinatoire des polyn\^omes 
orthogonaux g\'en\'eraux, Lecture Notes, UQAM, 1983.
\item{[V2]} G. Viennot, A combinatorial theory for general orthogonal
polynomials with extensions and applications, in Polyn\^omes 
Orthogonaux et Applications, Bar-le-Duc, 1984, Springer Lecture Notes 
in Mathematics, vol. 1171, 139-157.
\enddocument
\end